\documentclass[a4paper,12pt]{amsart}
\usepackage{amssymb,amsthm}
\usepackage[all]{xy}

\begin{document}

\newcommand{\opp}{\bowtie }
\newcommand{\po}{\text {\rm pos}}
\newcommand{\supp}{\text {\rm supp}}
\newcommand{\End}{\text {\rm End}}
\newcommand{\diag}{\text {\rm diag}}
\newcommand{\Lie}{\text {\rm Lie}}
\newcommand{\Ad}{\text {\rm Ad}}
\newcommand{\car}{\mathcal R}
\newcommand{\Tr}{\rm Tr}
\newcommand{\Spec}{\text{\rm Spec}}

\def\ge{\geqslant}
\def\le{\leqslant}
\def\a{\alpha}
\def\b{\beta}
\def\c{\chi}
\def\g{\gamma}
\def\G{\Gamma}
\def\d{\delta}
\def\D{\Delta}
\def\L{\Lambda}
\def\e{\epsilon}
\def\et{\eta}
\def\io{\iota}
\def\o{\omega}
\def\p{\pi}
\def\ph{\phi}
\def\ps{\psi}
\def\r{\rho}
\def\s{\sigma}
\def\t{\tau}
\def\th{\theta}
\def\k{\kappa}
\def\l{\lambda}
\def\z{\zeta}
\def\v{\vartheta}
\def\va{\varphi}
\def\x{\xi}
\def\i{^{-1}}

\def\mapright#1{\smash{\mathop{\longrightarrow}\limits^{#1}}}
\def\mapleft#1{\smash{\mathop{\longleftarrow}\limits^{#1}}}
\def\mapdown#1{\Big\downarrow\rlap{$\vcenter{\text{$\scriptstyle#1$}}$}}
\def\mapup#1{\Big\uparrow\rlap{$\vcenter{\text{$\scriptstyle#1$}}$}}

\def\ca{\mathcal A}
\def\cb{\mathcal B}
\def\cc{\mathcal C}
\def\cd{\mathcal D}
\def\ce{\mathcal E}
\def\cf{\mathcal F}
\def\cg{\mathcal G}
\def\ch{\mathcal H}
\def\ci{\mathcal I}
\def\cj{\mathcal J}
\def\ck{\mathcal K}
\def\cl{\mathcal L}
\def\cm{\mathcal M}
\def\cn{\mathcal N}
\def\co{\mathcal O}
\def\cp{\mathcal P}
\def\cq{\mathcal Q}
\def\car{\mathcal R}
\def\cs{\mathcal S}
\def\ct{\mathcal T}
\def\cu{\mathcal U}
\def\cv{\mathcal V}
\def\cw{\mathcal W}
\def\cz{\mathcal Z}
\def\cx{\mathcal X}
\def\cy{\mathcal Y}

\def\ccg{\mathfrak g}
\def\ccb{\mathfrak b}
\def\cch{\mathfrak h}
\def\cck{\mathfrak k}
\def\ccp{\mathfrak p}
\def\ccn{\mathfrak n}
\def\cct{\mathfrak t}
\def\ccl{\mathfrak l}
\def\cci{\mathfrak i}
\def\ccq{\mathfrak q}
\def\ccu{\mathfrak u}
\def\ccs{\mathfrak s}
\def\cco{\mathfrak o}
\def\ccd{\mathfrak d}
\def\cca{\mathfrak A}

\def\CR{\mathbb R}
\def\tz{\tilde Z}
\def\tl{\tilde L}
\def\tc{\tilde C}
\def\ta{\tilde A}
\def\tx{\tilde X}

\newtheorem{theorem}{Theorem}[section]
\newtheorem{lem}[theorem]{Lemma}
\newtheorem{cor}[theorem]{Corollary}
\newtheorem{prop}[theorem]{Proposition}
\newtheorem{thm}[theorem]{Theorem}
\newtheorem*{rmk}{Remark}
\newtheorem{ack}{Acknowledgement}
\newtheorem{eg}[theorem]{Example}
\newtheorem{defi}[theorem]{Definition}

\newenvironment{thmref}{\thmrefer}{}
\newcommand{\thmrefer}[1]{\renewcommand\thetheorem
{\protect\ref{#1}}\addtocounter{theorem}{-1}}

\author{Chuying Fang}
\address{Department of Mathematics, Massachusetts Institute of Technology, Cambridge, MA 02139, USA}
\email{cyfang@alum.mit.edu}
\title[Ad-nilpotent Ideals and Equivalence Relations]
{Ad-nilpotent Ideals and Equivalence Relations}

\begin{abstract}
In this paper  we study  ad-nilpotent ideals of a complex simple Lie algebra $\ccg$ and their connections with affine Weyl groups and nilpotent orbits. We  define a left equivalence relation for ad-nilpotent ideals based on their normalizer and generators,  and  prove that the equivalence relation is compatible with the left cell structure of affine Weyl group of $\ccg$ and  Lusztig's star operator for type $\tilde A_{n-1}$.
\end{abstract}

\maketitle

\section{Introduction}

Let $G$ be a quasi-simple Lie group over $\mathbb C$ and $\ccg$ be its Lie algebra. Fix a Borel subgroup $B$ of $G$ and let $\ccb$ be the Lie algebra of $B$ with nilradical $\ccn$.

A subspace of $\ccn$ is called an {\it ad-nilpotent ideal} if it is invariant under the adjoint action of $B$ (it is also called a $B$-stable ideal since it's stable under the adjoint action of the Borel subgroup $B$. ).  

Ad-nilpotent ideals have close relationship with affine Weyl groups  $\widehat W$ of $\ccg$. The connection between ad-nilpotent ideals and affine Weyl groups  was studied by  Cellini-Papi \cite{CP1} \cite{CP2},  Panyushev \cite{Pa1} \cite{Pa3}, Sommers \cite{So1}, Shi \cite{Sh1}, etc.    Cellini and Papi proved in \cite{CP1}\cite{CP2}  that there is a bijection between the ad-nilpotent ideals and certain elements of the affine Weyl group $\widehat W$ of $\ccg$. These elements are called minimal elements of $\widehat W$.    Shi showed in \cite{Sh2} that there is a natural bijection between  the set of ad-nilpotent ideals $\cca\ccd$  and  the dominant sign types,  and also the regions of the Catalan arrangement which are contained in the fundamental chamber.  Sommers studied in \cite{So1} both the ``minimal'' and  ``maximal'' elements of $\widehat W$ and showed that the set of maximal elements are in bijection with the set of strictly positive ideals of $\ccg$. Panyushev studied in \cite{So1} the ``maximal'' elements of $\widehat W$ and showed that they are in bijection with the set of strictly positive ideals. analyzed the generators and normalized of ad-nilpotent ideals in \cite{Pa1}\cite{Pa3} and obtained some statistics on them. 

Via restriction of the moment map on the cotangent bundle,  ad-nilpotent ideals also have natural connections with nilpotent orbits of $\ccg$ and related  representation theory. This involves the work of Kawanaka \cite{Kaw1}, Mizuno \cite{Mi}, Sommers \cite{So2} and Lawton \cite{Law}. It's an effective tool to analyze the structure theory of exceptional groups over finite field.   In addition, ad-nilpotent ideals are also used to study  some group cohomology  and its corresponding representation theory . In this paper,  our main interest is the equivalence relations on the affine Weyl groups and on the ideals.

For any ad-nilpotent ideal $I$, consider the map $G\times_B I \rightarrow \ccg$, which sends $(g, X) \in G \times I$ to $Ad(g)X$, where $Ad$ means the adjoint action of $G$. This is a restriction of the moment map and its image is the closure of one unique nilpotent orbit, which is denoted by $\co_I$. This orbit is  called the associated orbit of the ideal $I$.

There are various equivalence relations on the set of ideals $\cca\ccd$. In \cite{Mi}  Mizuno defined that  two ideals are equivalent if and only if they have the same associated orbit.  In \cite{So2}, Sommers defined an equivalence relation on $\cca\ccd$ based on some group cohomology.   We define a new equivalence relation $\sim_L$ on $\cca\ccd$ based on its generators and normalizers. This turns out to be related to the star operator of affine Weyl groups in the case of type $A$.

The left, right and two-sided cells of affine Weyl groups were introduced  by Kazhdan and Lusztig  \cite{KL} to study Hecke algebra and Kazhdan-Lusztig polynomials. Two elements lies in the same left cell are left equivalent. The star operation on $\widehat W$ was defined by Lusztig, which induces an $P_L$-equivalence relation on $\widehat W$. Indeed $P_L$-equivalence implies left equivalence.

For affine Weyl group of type $\tilde A_{n-1}$, Shi used  sign types to describe combinatorially the  cell structure of affine Weyl groups.  Shi showed in  \cite{Sh1} that the left cells (or two-sided cells) of $\widehat W$ are in bijection with the left cells (two-sided cells) of  sign types. Also it was proved by Lusztig  \cite{Lu} and Shi  \cite{Sh1} that the two-sided cells are parameterized by the partitions of $n$. 

In this paper, we study the left-equivalence relations of ad-nilpotent ideals and the equivalence relation of affine Weyl groups for type $\tilde A_{n-1}$. Our main result shows that for two left equivalent ad-nilpotent ideals, there exists two dominant elements of $\widehat W$ corresponding to these two ideals and they lie in the same left cell. On the other hand, if two elements of $\widehat W$ are $P_L$-equivalent, one of which is dominant, then the other one is also dominant and their associated ideals are left equivalent.  

The paper is organized as follows. We introduce the basic notations in section 2 and surveyed the known results about affine Weyl groups, the generator and normalizer of ad-nilpotent ideals. In section 3, we discussed  the moment maps and how to determine whether two ideals have the same associated orbits. We also introduced the left equivalence relations on ad-nilpotent ideals in this section. Section 4 mainly focused on dominant sign types.  In section 5, we state some commutative diagram about the ad-nilpotent ideals and prove our main theorem.

I would like to express my deep gratitude to my advisor David Vogan  for his guidance, warm encouragement and many useful suggestions. I would like to thank George Lusztig for his enjoyable lectures.  I  would also like to  thank Eric Sommers for both the email correspondence and conversations, pointing me the results of Shi and Lawton.

\section{Notations and Preliminaries}

Let $\mathbb N $ be the set of nonnegative integers and $\mathbb N^+$ be the set of positive integers.

Let  $\cch$  be a Cartan subalgebra of  $(\ccg, \ccb)$. We denote by  $\Delta$ the reduced root system associated to $(\ccg,\cch)$.   For each root $\alpha$, let $\ccg_{\alpha}$ be the corresponding root space in $\ccg$.    The Borel subalgebra $\ccb$ gives rise to  a positive root system $\Delta^+$ inside $\Delta$ and $\ccb$=$\cch \oplus  \bigoplus_{\alpha \in \Delta^+}\ccg_{\alpha} $.  Let $W$ be the Weyl group of $\ccg$.

Let $\Pi=\{\alpha_1, \alpha_2,\alpha_3, \dots, \alpha_n \}$ be the set of simple roots of $\Delta$. Then $Q= \oplus^n_{i=1} \mathbb Z \alpha_i $ is the root lattice. We set $V=\oplus_{i=1}^{n}\mathbb R \a_i =\cch^{\ast}_{\mathbb R}$.  The Killing form on $\ccg$ induces  a  $W$-invariant positive definite and symmetric  bilinear form on $V$, which is denoted by  $(\ , \ )$.  For each root $\alpha$,  $\alpha^\vee= \frac{2\alpha}{(\alpha,\a)} $ is the coroot for $\a$ and $Q^\vee= \oplus_{i=1}^p \mathbb Z \a_i^\vee$ is the coroot lattice.

Now we recall the definition of affine Weyl groups   based on  Kac's book  \cite{Kac}.

Let $\widehat{V}= V \oplus \mathbb R \delta  \oplus \mathbb R \lambda$. We extend the bilinear form  on $V$ to the bigger space $\widehat V $ by letting  $(\delta, \d)=(\delta, v)=(\lambda,v)=(\lambda,\lambda)=0$ for any $v \in V$ and $(\delta,\l)=1$. This is a nondegenerate symmetric bilinear form  on $\widehat V$ and by abuse of notation, we still denote it by $( \ , \ )$.

Let $\widehat \D= \{ \D + k \d \mid k \in \mathbb Z \}$  be the set of affine real roots. Let $$ \widehat \D^+ = (\D^++ \mathbb N \d) \cup (\D^- + \mathbb N^+ \d)$$ be the set of positive affine roots.  We denote by  $\a>0$ when  $\a $ is a positive affine root and $\a <0 $ when $\a \in \widehat \D$ is negative.

Let $\widehat \Pi=\Pi \cup \{-\a_0 +\d \}$ be the set of affine simple roots, where $\a_0$ is the highest root in $\D$.

For each $\a $ in $\widehat \D$, we define the reflection $s_\a$ on $\widehat V$ by $$s_\a(x)=x-\frac{2(\a,x)}{(\a,\a )}\a$$ for any $x \in \widehat V$. The affine Weyl group $\widehat W$ is generated by the set of  reflections $\{s_\a\}_{\a \in \widehat \D}$. For simplicity, we write $s_0$  for $s_{-\a_0 + \d}$ and $s_i$ for  the reflection $s_{\a_i}$. Let S be the set of simple reflections $\{s_0, s_1, \dots , s_n\}$.  Then $(\widehat W, S)$ is a Coxeter system.

Set \begin{gather*}\cc_0=\{\a \in V \mid (x,a)>0, \forall \a \in \Pi \text{ and } (x, \a_0)<1 \}, \\
\cc=\{ x\in V \mid (x,\a)>0, \forall \a \in \Pi \}.
\end{gather*}
We call $\cc_0$  the {\it fundamental alcove} and $\cc$ the {\it (open) fundamental chamber}.

For any  root $\a \in \D$ and any integer $k$, let  $H_{\alpha,k}=\{x \in V \mid (x,\a)=k\}$ be the hyperplane that is determined  by $\a$ and $k$.

Set \begin{align*} N(w) & = \{H_{\a,k} \mid H_{a,k} \space \text{ separates }  \cc_0 \text{ and } w(\cc_0) \} \\ &=\{ \a \in \widehat \D^+ \mid w \i (\a) \in -\widehat \D^+ \}.
\end{align*}

For each $w \in \widehat W $, let $l(w)$ be the cardinality of $N(w)$. We call $l$ the length function of $\widehat W$.





\begin{defi} We call a subset $\ci$ of $\D^+$ a (combinatorial) ideal of $\D^+$ if for any $\a \in \ci$,  $\b \in \D^+$ and $\a + \b \in \D$, we have $\a + \b \in \ci$.
\end{defi}

\begin{lem} The map $\ci \mapsto \oplus_{\alpha \in \ci} \ccg_{\a}$ gives a bijection from the set of ideals of $\D^+$ to the set of ad-nilpotent ideals of $\ccb$.
\end{lem}

\begin{rmk} Unless otherwise stated, we denote the ad-nilpotent ideal that corresponds to the ideal $\ci$ of $\D^+$ by $I$.\end{rmk}

Let $\ci$ be an ideal of $\D^+$. We call a root $\a$ a {\it generator} of $\ci$ if $\a \in \ci$ and for any positive root $\b \in \D^+$, $\a-\b \notin \ci$. The set of generators of $\ci$ (resp. $I$) is denoted by $\G(\ci)$ (resp. $\G(I)$).

For two ad-nilpotent ideals $I_1$ and $I_2$, the bracket relation is $$[I_1, I_2]=\{[X, Y] \mid \forall X \in I_1 \text{ and } \forall Y \in I_2 \}.$$

For any ideals $\ci_1$ and $\ci_2$ in $\D^+$, there is a similar a bracket relation: $$[\ci_1, \ci_2]=\{ \a+ \b \mid \a \in \ci_1, \b \in \ci_2 \text{ and } \a+\b \in \D \}.$$

Then $[\ci_1, \ci_2]$ is an ideal of $\D^+$ and corresponds to the ad-nilpotent ideal $[I_1, I_2]$.

Then we can derive inductively a sequence of ideals $\{ \ci^1, \ci^2, \dots , \ci^k \}$ from $\ci$ by letting $\ci^1= \ci$, $\ci^k=[\ci^{k-1}, \ci]$. This is a descending sequence of ideals.

Let $$\widehat W_{dom} = \{ w \in \widehat W \mid  w(\a) \geq 0 \text{ for all } \a \in \Pi \}.$$

Elements that lie in the subset $\widehat W_{dom}$ are called dominant. By the definition of $N(w)$,  the following conditions are equivalent:

(1) $w \in \widehat W $ is dominant.

(2) $N(w)$ is a subset of  $\cup_{k \geq 1 }(k \d - \D^+ )$.

(3) there's no hyperplane $H_{\a,0}$ separating $\cc_0$ from $w \i(\cc_0)$.

(4) $w \i (\cc_0)$ lies in the fundamental chamber of $V$.

The first two properties of $\widehat W$ come from the linear  action of $W$ on $\widehat W$ and the last two properties of $\widehat W$ come from the affine transformation of $\widehat W$ on $V$.

When $w$ is a dominant element of $\widehat W$, the set $ \{ \mu \in \D^+ \mid  \d - \mu \in N(w)\} $ is an ideal of $\D^+$. This induces  a map $\phi$ from $\widehat W_{dom}$ to the set of ideals in $\D^+$.  

In fact it's  a surjective map.  For any combinatorial  ideal $\ci$, there is a special dominant element $w$   that  corresponds to $\ci$. These elements were introduced in  \cite{CP1}.

\begin{defi} An element $w \in \widehat W$ is called minimal if

(1) $w(\a) \geq 0$ for all $\a \in \Pi$.

(2) If $\a \in \widehat{\Pi} $ and $w \i(\a)= k \d + \mu$, then $k \ge -1$.
\end{defi}

We denote the set of minimal elements by $\widehat W_{min}$. From part (1) of definition, it's obvious that $\widehat W_{min} \subset \widehat W_{dom}$. Moreover, we have the following result on the relations between the minimal elements and the ideals in $\D^+$.

\begin{prop} \cite[Prop2.12]{CP1}
There is a bijection between $\widehat W_{min}$ and the set of ideals in $\D^+$. The bijection is constructed as follows:

(a) For each $w \in  \widehat W_{min}$, its corresponding ideal is $\phi(w)$.

(b) For each ideal $\ci$ in $\D^+$,  the corresponding w is determined by the set of affine roots  $$N(w)=\cup _{k \ge 1}(k\d - \ci^k) \subset \widehat \D^+ .$$

\end{prop}

For each $w \in  \widehat W_{min}$, we denote by $\ci_w$ the ideal that is determined by the element $w$. Conversely, the minimal element that corresponds to $\ci$ is denoted by $w_{min} (\ci)$, or simply by  $w_\ci$. Because of the bijection between the set of ad-nilpotent ideals and the set combinatorial ideals, we have similar notations  $I_w$ and $w_{I}$ for the ad-nilpotent ideal $I$. We call the ad-nilpotent ideal $I_w$ the {\it first layer ideal} of $w$.

It is obvious that for any positive root  $\a$ in an ideal $\ci$, then $w_{\ci} \i(-\a+ \d) > 0$. 

Sommers gave a description of the set of generators of the ideal $\ci$ in  \cite[Cor6.3]{So1} and Panyushev showed independently in  \cite[Thm2.2]{Pa1}.

\begin{prop}\label{pan}
If $w \in  \widehat W_{min}$, then any positive root $\a \in \ci$ is a generator of the ideal $\ci_w$ if and only if $w(\a -\d)$ is an affine simple root in $\widehat{\Pi}$.
\end{prop}

Given an ad-nilpotent ideal $I$, since $\ccb$ normalize $I$, the normalizer of $I$ must be a parabolic subalgebra containing $\ccb$. Therefore the normalizer is determined by the simple root Levi subalgebras.

Suppose $\ccl(\a)$ is the Levi subalgebra in $\ccg$ corresponding to a simple root $\a \in \Pi $. Namely $\ccl(\a)= \cch \oplus \ccg_{\a} \oplus \ccg_{-\a}$. Panyushev showed in  \cite[Thm2.8]{Pa3} that:

\begin{prop}\label{pa}
If $I$ is an ad-nilpotent ideal of $\ccb$, then $\ccl(\a) \subset N_{\ccg}(I)$ if and only if $w_I(\a) \in \widehat{\Pi}$, where $N_{\ccg}(I)$ denotes the normalizer of $I$ in $\ccg$.
\end{prop}

\section{The Moment Maps}\label{mom}

Let $\cb$ be the set of all Borel subalgebras in $\ccg$, $T^{*}\cb$ the cotangent bundle over $\cb$ and $\cn$ the nilpotent cone of $\ccg$.  Following  \cite[Lem1.4.9]{CG},  if one identifies the Lie algebra $\ccg$ with its dual $\ccg^*$ via the Killing form of $\ccg$, the moment  map $$ m: T^*(\cb) \longrightarrow \cn $$
is equivalent to
$$ m: G \times_B \ccn \longrightarrow \cn. $$
Therefore, for each ad-nilpotent ideal $\ci \subset \ccn$, $G \times_\cb I $ may be  considered as a $G$-equivariant subbundle of cotangent  bundle $T^* \cb$. The image   of $G \times _B I$ under the moment map is the closure of a nilpotent  orbit in
$\cn$.

Let $N(\ccg)$ be the set of nilpotent orbits of $\ccg$. The moment map $m$ induces a map  $$p: \cca\ccd \rightarrow N(\ccg)$$ by
sending $I \in \cca\ccd$ to a nilpotent orbit  $\co$, where the image of $G \times_B I$ under the moment map  is the closure of $\co$. The nilpotent orbit $\co$ is called the associated orbit of the ideal $I$.

It is easy to see that the map is surjective. Indeed, if $e \in \co$, by Jacobson Morozov theorem, there exists an $\ccs\ccl_2$-triple $\{e, h, f \}$ with $h \in \cch$. We conjugate the triple with elements of $G$ such that $\a(h) \ge 0  $ for positive root  $a$. Then $I_{h}:=\oplus_{i \ge 2} \ccg_{i}$, where $\ccg_{i}=\{X \in \ccg \mid [h, X]=iX \}$   is an ad-nilpotent ideal and    $\co_{I_h}=\co$.

Concerning this map $p$, there  are several   natural questions to ask:

When two different ideals give rise to the same nilpotent orbit?

Is it possible to describe combinatorially  the equivalence relation on minimal elements of the affine Weyl group corresponding to $p(I_{w_1})=p(I_{w_2})$?

In the special case that the two ideals differ by a single positive root, then the following  result from  \cite{GS} that partially answers the first question.

\begin{prop}\label{sommmer}
Let $I$ be an ad-nilpotent ideal of $\ccn$. Suppose that $I$ is stable under the adjoint action of $\ccl(\a)$ for some simple root $\a \in \Pi$ and $\b$ is a generator of $I$ such that $s_{\a}(\b)> \b $.  Let $J$ be the ad-nilpotent ideal such that $I=J\oplus \ccg_\b$. Then  $I$ and $J$ has the same image under $p$.
\end{prop}

If two ideals $I$ and $J$ satisfies the condition of Proposition \ref{sommmer},  we denote  by $I \sim_L  J$ or $ J \sim_L I$.

\begin{defi}\label{left}
Two ad-nilpotent ideals $I$ and $J$ are called left equivalent if either $I=J$ or there exists a sequence of ideals $I_1= I, I_2, \dots, I_k =J$, such that $I_i \sim_L I_{i+1}$ for $i= 1, \dots, k-1$.
\end{defi}

Proposition \ref{sommmer} is a sufficient condition for  two  ad-nilpotent ideals to have  the  same associated orbit. Since the moment map sends the left-equivalent classes of ideals to the same orbit, one may ask whether it is true that two ad-nilpotent ideals  have the same associated orbits only if these two ideals are left equivalent.  For some cases, we can give an affirmative answer.

\begin{eg}
When $\ccg$ is of type $G_2$. Suppose that $\a_1$ is the short simple root, and $\a_6$ the long root. The other four roots are $\a_2= \a_1+ \a_6$, $\a_3= 2\a_1+ \a_6$, $\a_4= 3\a_1+ \a_6$, $\a_5= 3 \a_1+ 2\a_6$. There are eight ideals and five equivalence classes:

1. The zero ideal, which corresponds to the zero orbit;

2. The ideal $\ccg_{\a_3}\oplus \ccg_{\a_4}\oplus \ccg_{\a_5}$, which corresponds to the 8-dimensional orbit;

3. The maximal ideal, which corresponds to the principal orbit;

4. $ \ccg_{\a_5}  \sim_L  \ccg_{\a_5}\oplus \ccg_{ \a_4}$, which corresponds to the minimal orbit;

5. $\ccg_{\a_1} \oplus \ccg_{\a_2}\oplus \ccg_{\a_3}\oplus \ccg_{\a_4}\oplus \ccg_{\a_5} \sim_L \ccg_{\a_2}\oplus \ccg_{\a_3} \oplus \ccg_{\a_4} \oplus \ccg_{\a_5} \oplus \ccg_{\a_6} \sim_L \ccg_{\a_2}\oplus \ccg_{\a_3} \oplus \ccg_{\a_4}\oplus \ccg_{\a_5}$, which corresponds to the subregular orbit.
\end{eg}

In this example above, the left equivalence relation on the ideals completely determines the image of the moment map. For classical cases other than  type $A$, Proposition \ref{sommmer} is not necessary to determine the fiber of the moment map. The simplest counter example is in type $B_3$.
\begin{eg}
Let $\ccg$ be $\ccs\cco(7)$. Suppose that $\a_1$, $\a_2$ are the long simple roots of $\ccg$ and $\a_3$ is the short simple root. Let $I$ be the ad-nilpotent ideal with generators  $\a_2+\a_3, \a_1+\a_2$ and $J$ be the ad-nilpotent ideal with generators by $\a_2+2\a_3$ and $\a_1+\a_2$. Then $I= J\oplus \ccg_{\a_2+\a_3}$, and $I, J$ have the same associated orbit. But no Levi subalgebra of $\ccg$ normalizes $I$ except $\cch$. $I$ and $J$ don't satisfy the condition of Proposition \ref{sommmer}.
\end{eg}


In the type $A$ case, we haven't found such a counterexample yet. We would prove some results in the following section which may shed some light on whether Proposition \ref{sommmer} completely determines the fiber of the moment map.

Before directly approaching this problem, let us   make a short digression and turn to the geometric description of the ideals via sign types. It turns out sign types are very useful for us to  understand affine Weyl groups and   are also  connected with   ad-nilpotent ideals.

\section{Sign types}
Recall that  for any positive root $\a \in \D^+$, $H_{\a, k}$ is the hyperplane defined by  $H_{\a, k}= \{v \in V \mid (v,\a)=k  \}$.   We denote three regions that are separated  by two hyperplanes $H_{\a, 0}$ and $H_{\a, 1}$ by: \begin{gather*}
\tilde H_{\a, +}= \{ v \in V \mid (v,\a )>1\}, \\
\tilde H_{\a, 0}= \{ v \in V \mid  0<(v,\a )<1 \}, \\
\tilde H_{\a, -}= \{ v \in V \mid (v,\a )<0\}.
\end{gather*}

The non-empty connected simplex of $ V- \cup _{\a \in \D^+, k=0,1} H_{\a, k}$  is called  a {\it sign type} of $V$.  The definition of sign type was first introduced by Shi in  \cite{Sh1} to study the left cells and two-sided cells of affine Weyl group of type $\tilde A_{n-1}$ and later was generalized by him to other classical groups in  \cite{Sh2}. Notice that in Shi's original definition, he used the coroot system  of $\ccg$ while here we use the root system.

We denote  the set of all sign types   by $\cs$. For any sign type $s \in \cs$, it lies in one of the three regions $\tilde H_{\a, \e}$, where $\e \in \{0, +, -\}$. Therefore each sigh type has the form $s= \cap_{\a \in \D^+, \e_{\a} \in \{+, 0, -\}}  H_{\a , \e_{\a}}$, with $\a \mapsto \e_\a$  mapping  $\D^+$ to the index set $\{0, +, - \}$. This  map completely determines the sign type.


Since the walls $H_{\a, 0 }$ for $\a \in \D^+ $ are used to define  sign types, each sign type either lies in the fundamental chamber or has no intersection with  the fundamental chamber.  We call the special sign types that lie inside the fundamental chamber $\cc$ {\it dominant sign types}.  The set of dominant sign types in $\cs$ is denoted by $\cs_{dom}$.  It's clear from the definition that the sign type $s$ is dominant if and only if it lies in the region $\tilde H_{\a, \e}$, for each $\a  \in \D^+$ and $\e \in \{
+, 0 \}$.

\begin{prop}\label{shi} \cite{Sh3}
The map $$I \mapsto R_I:=\{ x \in \cc \mid (x,\a)>1 \text{ if } \ccg_\a \subset I \text{ and } 0<(x, \a)<1 \text{ if } \ccg_\a \nsubseteq I \}$$ gives a bijection from the set $\cca \ccd$ of  ad-nilpotent  ideals in $\ccb$ to  $\cs_{dom}$.
\end{prop}

Sign types are also  closely related to affine Weyl groups. For each $w \in \widehat W$,  $w \i $ maps the fundamental alcove $\cc_0$ to another alcove.  Thus $w \i (\cc_0)$ is contained in a unique sign type  $s$. We obtain a map $\widehat W \rightarrow \cs$ by sending  $w$ to  $s$. It turns out that dominant elements in $\widehat W$ are mapped to dominant sign types. Namely, when restricted to $\widehat W_{dom}$, we have $\widehat W_{dom} \rightarrow \cs_{dom}$.

Recall that for each dominant element $w \in \widehat W$, one can associate an ad-nilpotent ideal $I= \{ \a \in \D^+ \mid w(\d -\a ) \leq 0 \}$. This is equivalent to say that the hyperplane $H_{\a, 1}$ separates $w \i (\cc_0)$ from $\cc_0$ if and only if the positive root $\a$ lies  in the ideal $I$. In other words, $w \i(\cc_0)$ lies in the region $\tilde H_{\a, +}$ for any positive root $\a \in I$  and  lies in the region $\tilde H_{\a, 0}$ for any positive root $\a$ that is not in the ideal $I$.

Comparing this with the map $\widehat W \rightarrow \cs$ defined in the preceding paragraph and the bijection between dominant sign types and ad-nilpotent ideals in Proposition \ref{shi},  we have a commutative diagram.
\[ \tag{1.4.1} \xymatrix{
\widehat W_{dom} \ar[r] \ar[d] &  \cca \ccd \ar[ld]^{\simeq} \\
\cs_{dom} &}
\]

Under this commutative diagram, it's possible for us to explore the  relation between the equivalent classes of affine Weyl group and the equivalent classes of  ad-nilpotent ideals given by Definition \ref{left} in further detail.

\section{The Case of Type $\tilde A_{n-1}$}

Now we focus on the affine Weyl group of  type $\tilde  A_{n-1}$. In this section, $\widehat W$ denotes the  affine Weyl group of type $\tilde A_{n-1}$.

For $W$,  there's a combinatorial  description of the affine Weyl group elements.
We regard $\widehat W$ as a set of permutations on $\mathbb Z$ as follows  \begin{align*} \widehat W=\{\s: \mathbb Z \rightarrow \mathbb Z \mid \text{ }& \s(i+n)=\s(i) \text{
for }i \in \mathbb Z, \\ & \sum_{t=1}^n (\s(t)-t)=0(\text{mod }  n) \}.
\end{align*}

The simple reflection $s_0, s_1, \dots s_{n-1}$ can be taken to the permutations as follows. For $0 \le j \le n-1$, $j$ is mapped to:

\[
s_i(j)=\begin{cases}
    j+1,  & \text{if } j \equiv i(\text{mod } n); \\
    j-1,& \text{if } j \equiv i+1(\text{mod } n);  \\
    j,   & \text{  otherwise }.
    \end{cases}
\]

We denote by $<$ the Bruhat order on $\widehat W$ and  denote the set of partitions of $n$ by $P(n)$.

Let $[n]$ be the index set $\{1, 2, \dots, n\}$. For any $w \in \widehat W$,   the partial order  $\succ_w$  of $[n]$ is defined by $i \succ_w j$ if and only if either  $ i< j$ and $w(i)>w(j)$ or $i> j$ and $w(i)>w(j)+n$. A chain of $[n]$ is a sequence of integers $\{i_1, \dots, i_h\}$ with the order $i_1 \succ_w i_2 \dots \succ_w i_h$. A $k$-family is a subset of $[n]$ such that there's no chain of length $k+1$ in this subset.

Let $d_k$ be the maximal cardinality of a $k$-family of $[n]$. Then  $d_1 \le d_2 \le \dots \le d_n=n$. Let $\l_1=d_1$, $\l_j=d_j-d_{j-1}$ for $1 \le j\le n$. By a theorem of Green  \cite{Gr}, $\l_1 \ge \l_2 \dots \ge \l_n$ and $\sum_{i=1}^n\l_i=d_n=n$, which gives us a partition  of $n$. Let $\mu=\{\mu_1 \ge \mu_2 \dots \ge \mu_n \}$ be the conjugate partition of $\l=\{\l_1\ge \dots \ge \l_n \}$, meaning that for any $h=1, \dots, n$, $\mu_h$ is equal to the number of parts in $\l$ of size $h$.  This gives us a map $\phi: \widehat W \rightarrow P(n)$ with  $\phi(w)= \mu$. This map is related to two-sided cells of $\widehat W$ (The reason that we choose the conjugate partition $\mu$ instead of $\l$ is  to get  a commutative diagram 1.5.1  below).


Set \begin{gather*} L(w)= \{s_j \in S \mid s_j w < w \}=\{\a_j \in \widehat{\Pi} \mid w \i(\a_j) <0\}, \\ 
\end{gather*}

We follow the usual definition of left cells, right cells and two-sided cells in  \cite{KL}.

For any two elements $w, u \in \widehat W$, we denote by $w \sim_L u$ (resp, $w \sim_R u $; or $w \sim_{LR} u$) if $w$ and $u$ lie in the same left cell (resp, right cell; or two sided cell) of $\widehat W$.

It was proved by Lusztig in  \cite{Lu} and Shi in  \cite{Sh1} that the  two-sided cells of $\widehat W$ is parameterized by $P(n)$ and these two-sided cells coincide with fibers of $\phi$. In particular, to describe the fiber of the map, Shi  \cite{Sh1}  constructed a map $\Phi: \cs \rightarrow P(n)$ such that   the following is a commutative diagram:
\[
\xymatrix{ \widehat W \ar[r] \ar[d]_-{\phi} & \cs \ar[ld]^-{\Phi}\\
P(n)}
\]


For the map $\widehat W \rightarrow \cs$, Shi proved in  \cite[Chap 18]{Sh1} that if two affine Weyl group elements are mapped to the same sign type, then these two elements lie in the same left cell. Because of the commutative diagram 1.4.1, two dominant elements of $\widehat W$ with the same associated ideal also lie in the same left cell.

When $\ccg=\ccs\ccl(n)$,  the set of nilpotent orbits of $\ccg$ is parameterized by $P(n)$ (see  \cite{CM}) and  the  moment map of $\ccg$  gives us   $p: \cca\ccd \rightarrow P(n)$.  Since the set of dominant sign types  $\cs_{dom}$ is in bijection with the set of ad-nilpotent ideals,  there might be some  relation between  the map $\Phi|_{\cs_{dom}}$ and the  map $p$.

Indeed the work of Shi (see  \cite{Sh1}) and the work of Lawton (see  \cite{Law}) imply that there exists a commutative diagram in the following form:
\[
\tag{1.5.1}
\xymatrix{{\widehat W}_{dom}\ar[r] \ar[d]^-{id} & \cs_{dom} \ar[r]^-{\Phi|_{\cs_{dom}}} \ar[d]& P(n) \ar[d]^{id}\\
\widehat W_{dom} \ar[r] & \cca \ccd \ar[r]^{p} & P(n)}
\]

Recall that  Proposition \ref{sommmer} gives us a criterion to determine when two ad-nilpotent ideals have the same associated orbits. One might hope that this criterion is related to  the cell structure of affine Weyl groups and sign types. Indeed, it's  related to the left star operation on $\widehat W$. 

Let $s, t$ be two simple reflections of $\widehat W$ such that $st$ has order 3.

\begin{defi}
\begin{gather*} D_L(s,t)= \{w \in \widehat W \mid L(w) \cap \{s, t\} \text { contains only one element }  \}, \\
\end{gather*}
\end{defi}  If $w$ is in the set $D_L(s,t)$, then $\{tw, sw\} \cap D_L(s,t )$ contains only one element, which is denoted by ${}^*w$. Then the map that sends $w$ to ${}^*w$ defines an involution on $D_L(s,t)$ and is called a left star operation.  This involution depends on the two simple reflections $s$ and $t$. For example, if $w$ lies both in the sets $D_L(s_1, t_1)$ and $D_L(s_2, t_2)$ for simple reflections $s_1$, $s_2$, $t_1$, $t_2$, then the involution $*$ may give two  different ${}^*w$.

The left star operation generates  another equivalence relation  on $\widehat W$. The element $w$ is $P_L$ equivalent to $w'$ if and only if there is a succession of affine Weyl group elements $w_1=w, w_2, \dots, w_n =w'$ such that each $w_i$ lies in some set $D_L(s_i, t_i)$ and $w_{i+1}={}^*w_i$ under the left star operation  in $D_L(s_i, t_i)$. We denote this equivalence relation by $\sim_{P_L}$.

It is proved in  \cite{KL} that $w$ and ${}^*w$ lie in the same left cell, therefore  $w \sim_{P_L} w'$ implies that  $w \sim_L  w'$.

Now we can come to the proof of the main theorems.

\begin{theorem}\label{pl}
Suppose $I_1 $ and $I_2 $ are two ideals and $I_1 \sim_L I_2$, then there exists two  dominant elements $w_1$ and $w_2$, such that $I_{w_1}=I_1$, $ I_{w_2} =I_2$ and $w_1 \sim_L w_2$.
\end{theorem}

Proof. Since any  elements of $\widehat W_{dom}$ with the same  associated ideal lie in the same left cell, by the definition of left equivalence relation  for the  ideals, it suffices to prove the statement when $I_1$ and $I_2$ differs by only  one positive root.

Suppose that $I_2 \subset I_1$, $\a$ is a simple root that normalize $I_1$ and $\b$ is a generator of $I_1$, s.t. $I_1= I_2 \oplus \ccg_\b$ and $s_\a(\b)> \b$.  Let $w_{I_1}$ be the minimal element in $\widehat W$ that corresponds to $I_1$ (By Proposition \ref{shi}). Recall from Proposition \ref{pa}, $w_{I_1}$ maps $\a$ to a simple root in $\widehat \Pi$, which we denote by $\a_i$. On the other hand, since $\b $ is a generator of $I_1$, by Proposition \ref{pan}, $w_{I_1}(\b -\d ) =\a_j $ for some simple root $\a_{j}$. The condition that $s_\a(\b) > \b$ implies that $(\a, \b) <0$. We denote $s_{\a_i}$ by $s_i$ and $s_{\a_j}$ by $s_j$. Then $s_i, s_j$ are two simple reflections of $\widehat W$. Since the inner product on $\widehat V$ is invariant under $\widehat W$, $(\a_i, \a_j)<0$  and $s_is_j$ has order 3. The fact that $w_{I_1}(\a)=\a_i$ implies that $s_i w_{I_1} > w_{I_1}$ and $w_{I_1}(\b -\d)= \a_j$ implies that $s_j w_{I_1}< w_{I_1}$.  Therefore, $w_{I_1}$ is an element that lies in the set $D_L(s_i, s_j)$. In this case, we can determine ${}^*w_{I_1}$ explicitly.  Indeed it's clear  that $s_j(s_jw_{I_1})=w_{I_1} > s_j w_{I_1}$ and $s_i(s_jw_{I_1}) <s_jw_{I_1}$ because:
$$(s_jw_{I_1}) \i (\a_i)= w \i s_j(\a_i)= w_{I_1} \i (\a_i +\a_j)= \a +\b -\d <0$$
Hence  ${}^*w_{I_1}= s_jw_{I_1}$. In addition, $l(s_jw_{I_1})< l(w_{I_1})$ and $N(w_{I_1})=N(s_jw_{I_1}) \cup (\b-\d)$. Hence $N(s_jw_{I_1})$ is a subset of $\cup_{k \ge 1} (k\d - \D^+)$, which implies that $ s_jw_{I_1}$ is still dominant. We have $I_{s_jw_{I_1}} \oplus \ccg_\b= I_1$, i.e. $I_{s_jw_{I_1}}=I_2$. \qed

If two dominant elements in $\widehat W$ are $P_L$-equivalent, we can prove a converse version of Theorem \ref{pl}.   First we need to prove a lemma that states a special property for the $P_L$ equivalence classes of elements in $\widehat W_{dom}$.

\begin{lem}\label{gen}
Suppose that $w$ is dominant and $\a$ is a simple root in $\widehat \Pi$. If $w\i(\a)<0$, then $w\i(\a)= -k \d + \b$, where $k $ is positive and $\b$ is an element  of $I_w$. In particular, if $k=1$, then $\b$ is a generator of $I_w$.
\end{lem}

Proof. The element $w$ is dominant implies that $w \i(\a)= -k \d+ \b$, where $k  \ge 1$ and $\b \in \D^+$.  Since $w(\d- \b)= -\a - (k-1) \d <0$, $\b$ lies in the first layer ideal $I_w$.
If $k=1$, for any $\g \in \D^+$, $w(\d -(\b-\g))=w(\d-\b) + w(\g)= -\a+ w(\g) $. Since $w(\g)$ is a positive affine root and $\a$ is simple, $w(\d -(\b-\g))>0$ and $\b -\g$ does not belong to the ideal $I_w$, which means that $\b$ is indeed a generator of $I_w$.  \qed

\begin{lem}\label{gen1}
If $w \in \widehat W_{dom}$ and $w$ lies in the set $D_L(s_i, s_j)$ for some simple reflections $s_i $ and $s_j$ with $s_is_js_i= s_js_is_j$, then ${}^*w$ is also a dominant element.
\end{lem}

Proof. Any element $u \in \widehat W$ is dominant if and only if $N(u) $ is contained in $\cup_{k \ge 1} (k\d - \D^+)$. There are two possibilities for ${}^*w$. The first case is  $l({}^*w)= l(w)-1$.  Then
$N({}^*w)$ is a subset of $N(w)$ and the fact that $w$  is dominant implies  ${}^*w  $ is also dominant.

The other case is $l({}^*w)=l(w)+1 $. In this case, by the symmetry of  $i$ and $j$, suppose ${}^*w =s_jw$.  Then $s_jw >w, s_iw<w$ and $s_is_jw >s_jw$. Let $a_i$, $a_j$ be the two simple roots of $\widehat \Pi$ that correspond to $s_i$ and $s_j$.  By the dominance property of $w$ and the fact that $s_iw <w$,  we have  $w \i(a_i)= t \d -\b$, where $t$ is a positive integer and $\b \in \D^+$.
Similarly $s_jw >w $ means that  $w\i (\a_j)>0$, namely $w \i(\a_j)= k\d + \g$, where $k \ge 0$ and $\g \in \D^+$ or $k \ge 1$ and $\g \in \D^-$.

Case (a): $k=0$ and $\g \in \D^+ $.
From the facts that $s_is_j w> s_jw$ and $\a_j$ has the same length as $\a_i$, $w \i(\a_i+\a_j)= (s_jw) \i(\a_i)>0$. On the other hand, $w\i(\a_i+\a_j)= \b+\g -k\d <0 $, which is a contradiction.

Case (b) : $k>0$ and $\g \in \D^+$.  Since $w$ is dominant, $w(\g)$ is a positive root in $\widehat \D^+$. Also $w$ fixes $k\d$. Therefore, $\a_j=w(k\d +\g) $ can not be a simple root in $\widehat \D^+$. Contradiction.

The only possible form for $w \i(\a_j) $ is $k\d +\g$, where $k \ge 1$ and $\g \in \D^-$. Notice that $N({}^*w)$ is the union of $N(w)$ and $w \i(\a_j)$, then ${}^*w$ is again dominant. \qed

\begin{lem} \cite[Thm3.5]{Pa3}\label{gen2}
Suppose $w\in \widehat W$ is dominant, and let $I_w$ be the corresponding ideal of $w$. If there exist a simple root $\a \in \Pi$, such that $w(\a) \in \widehat \Pi$, then $\ccl(\a)$ normalize $I_w$.

\end{lem}

\begin{theorem}
Suppose that $w_1$  is  dominant  and $w_1 \sim _{P_L} w_2$. Then $w_2$ is also dominant and the  corresponding ideals of $w_1$ and $w_2$ are left equivalent, i.e. $I_{w_1} \sim _{L} I_{w_2}$.
\end{theorem}

Proof. The dominance of $w_2$ is shown in Lemma \ref{gen1}. By definition of $P_L$-equivalence relation, the problem can be reduced to the case when $w_1 $ lies in the set $D_L(s_i, s_j)$ for some $i, j \in \{0,1, \dots n\} $ and $w_2= {}^*w_1$. In addition, by the symmetry of $w_1$ and $w_2$, as well as the symmetry of $i$ and $j$, we may assume that $l({}^*w_1)<l(w_1)$ and ${}^*w_1= s_i w_1$. Then $s_iw_1< w_1, s_jw_1>w_1$ and $s_js_iw_1<s_iw_1$. This means that $w_1 \i (\a_i)<0, w_1\i \a_j>0$, and $w_1^{-1}(\a_i+\a_j)= (s_iw_1) \i (\a_j)<0$. By lemma \ref{gen}, $w_1\i(\a_i)= \b - k\d $, where $\b$ is an element of $I_{w_1}$, and $k\ge 1$.

If $k>1$, then $k\d - \b$ belongs to $N(w_1)$ and any $t \d - \b $ for $1 \le t <k$ also belongs to $N(w_1)$. That means $k\d -\b$ is the only element that lies in  $N(w_1)$, but not in $N(s_iw_i)$.  Since the definition of the corresponding ideals of $s_iw_1$ and $w_1$  involves only the first layer of positive roots in $s_iw_1$ and $w_1$, $s_iw_1$ and $w_1$ have the same first layer ideal, i.e. $I_{w_1}= I_{s_iw_1}$.

Suppose $k=1$. In this case, from Lemma \ref{gen}, $\b $ is a generator of $I_{w_1}$. In addition, we have known that  $w_1\i (\a_j)>0$ and $w_1\i(\a_i+\a_j)<0$. There are four possible cases for $w_1\i (\a_j)>0$.

Case(a): $w_1\i (\a_j)= k\d + \b_1$, where $k \ge 1$ and $\b_1 \in \D^+$.

In this case, $w_1\i (\a_i+\a_j)= (k-1)\d + \b+\b_1 >0$. Contradiction.

Case(b): $w_1\i (\a_j)=k \d - \b_1$, where $k \ge 2$ and $\b_1 \in \D^+$.

Similar to case (a), $w_1\i (\a_i+\a_j)= (k-1)\d + \b- \b_1>0$. Contradiction.

Case(c):  $w_1\i (\a_j)= \d-\b_1$, where $\b_1 \in \D^+$.

In this case, $w_1\i (\a_i+\a_j)=\b- \b_1 \in \D^-$. Suppose that $\b-\b_1 = -\g$, where $\g $ lies in $\D^+$. Then $\b_1=\b+ \g$, which implies that $\b_1$ also lies in ideal $I_{w_1}$ since $\b$ is a generator of $I_{w_1}$.  This contradicts the fact that $w_1(\d-\b_1)=\a_j$, which is positive.

Case(d): $w_1\i (\a_j)= \b_1$, where $\b_1 \in \D^+$.

This is the only possible form for $w_1\i (\a_j)$. In this case, $\b_1$ is a simple root of $\D^+$. Indeed, if $\b_1$ is not simple, say $\b_1= \g_1 + \g_2$, where $\g_1$ and $\g_2$ are two positive roots in $\D^+$, then $w_1(\b_1)= w_1(\g_1)+ w_2(\g_2)$, where $w_1(\g_1)$ and $w_2(\g_2)$ are two positive roots in $\widehat \D^+$. This contradicts the fact that $\a_j=w_1(\b_1)$ is a simple affine root. From  lemma \ref{gen2}, $\ccl(\b_1)$ normalize the ideal $I_{w_1}$. Also $\b$ is a generator of $I_{w_1}$ and $(\b, \b_1)= (\a_i,\a_j)<0$. Since $N(w_1)$ is equal to $N(s_iw_1) \cup (\d- \b)$, $\b$ appears in the first layer ideal of $w_1$, but not the first layer ideal of $s_iw_2$.  By the definition of left equivalence of ideals, $I_w$ is left equivalent to $I_{s_iw}$.  \qed

\

We have shown above that if two dominant  elements  in $\widehat W$ are $P_L$-equivalent, then their corresponding ideals are left equivalent (see Definition \ref{left}) and also two left equivalent ideals would give us two $P_L$ equivalent elements. It was proved in  \cite[Chap 19]{Sh1} that two dominant sign types lie in the same two-sided cell if and only if they lie in the left cell. Therefore, if we can prove that dominant elements in the same left cell give rise to left equivalent ideals, then Proposition \ref{sommmer} is the necessary condition to determine the  two-sided cell structure of sign types. However, Shi defined two operations for elements in the same left cell. The first one is related to $P_L$-equivalence relation, which we showed above that it does give two left equivalent ideals. The other one is called raising layer  operation, which we can not get its relation with the left equivalence relation of ad-nilpotent ideals. 

\bibliographystyle{amsalpha}

\end{document}